\documentclass[12pt]{amsart}
\usepackage{latexsym}
\usepackage[all]{xy}
\usepackage{amsfonts}
\usepackage{amsthm}
\usepackage{amsmath}
\usepackage{amssymb}

\def\<{{\langle}}
\def\>{{\rangle}}

\def\note#1{{}}

\def\note#1{}

\def\ca{{\mathcal A}}

\def\beq{\begin{equation}}
\def\eeq{\end{equation}}

\def\ot{{\otimes}}

\newcounter{zlist}

\def\Label#1{\label{#1}\ifmmode\llap{[#1] }\else
\marginpar{\smash{\hbox{\tiny [#1]}}}\fi}
\def\Label{\label}

\newtheorem{proposition}{Proposition}[section]

\newtheorem{theorem}[proposition]{Theorem}
\theoremstyle{definition}

\theoremstyle{remark}

\newcounter{c}

\newcommand{\etyk}[1]{\vspace{-7.4mm}$$\begin{equation}\Label{#1}
\addtocounter{c}{1}}
\renewcommand{\]}{\ifnum \value{c}=1 $$\else \end{equation}\fi}
\begin{document}

\title{New constructions of Yang--Baxter systems}
\author{Florin F. Nichita}
\author{Deepak Parashar}
\address{Institute of Mathematics "Simion Stoilow" of the Romanian Academy, P.O. Box 1-764, RO-014700 Bucharest, Romania}
\address{American University of Kuwait, The College of Arts and Science, Safat 13034, Kuwait}
\email{Florin.Nichita@imar.ro}
\address{Mathematics Institute, University of Warwick, Coventry CV4 7AL, U.K.}
\address{Max-Planck-Institut f\"ur Mathematik, Vivatsgasse 7, 53111 Bonn, Germany}
\email{D.Parashar@warwick.ac.uk, deepak@mpim-bonn.mpg.de}
\begin{abstract}
The quantum Yang--Baxter equation admits generalisations to systems
of Yang--Baxter type equations called Yang--Baxter systems. Starting
from algebra structures, we propose new constructions of some
constant as well as the spectral-parameter dependent Yang--Baxter
systems. Besides, we also present explicitly the commutation algebra
structure associated to the constant type in dimension two.

\noindent \it{AMS Contemporary Math. 442 (2007) 193--200.}
\end{abstract}
\maketitle

\section{Preliminaries}

{\em Yang--Baxter systems} emerged from the study of quantum
integrable systems, as generalisations of the quantum Yang--Baxter
equation (QYBE) related to nonultralocal models \cite{h1,hk}. In
deriving the relations for the quantum monodromy matrices, it is
common to assume that the quantum integrable models under
investigation are ultralocal, i.e. the quantised Lax operators
corresponding to different sites of the lattice commute. This is no
longer the case for nonultralocal models where the relations for the
monodromy matrices satisfy certain conditions involving a collection
of several Yang--Baxter type equations which form the so-called
Yang--Baxter systems.

It is convenient to describe Yang--Baxter systems in terms of
Yang--Baxter commutators. Let $V$, $V'$, $V''$ be finite dimensional
vector spaces over a field $k$. Consider three linear maps $R: V\ot
V' \rightarrow V\ot V'$, $S: V\ot V'' \rightarrow V\ot V''$ and $T:
V'\ot V'' \rightarrow V'\ot V''$. Then, a {\em constant Yang--Baxter
commutator} is a map $[R,S,T]: V\ot V'\ot V'' \rightarrow V\ot V'\ot
V''$ defined by \beq [R,S,T]:= R_{12} S_{13} T_{23} - T_{23} S_{13}
R_{12} \eeq where $R,S,T$ are understood as constant $n^2 \times
n^2$ matrices ($n$ being the dimension of V) and $[R,S,T]$ a $n^3
\times n^3$ matrix. Note that $[R,R,R] = 0$ is just a short-hand
notation for writing the {\em constant} QYBE \beq \label{qybe}
R_{12} R_{13} R_{23} = R_{23}R_{13} R_{12}\eeq Similarly, a {\em
coloured} or spectral-parameter dependent {\em Yang--Baxter
commutator} is the $n^3 \times n^3$ matrix
\begin{eqnarray}
[[R,S,T]] = && [[R,S,T]](u,v,w):= \nonumber\\
&& R_{12}(u,v) S_{13}(u,w) T_{23}(v,w) - T_{23}(v,w) S_{13}(u,w)
R_{12}(u,v)
\end{eqnarray} where the $n^2 \times n^2$
matrices $R,S,T$ now depend upon the spectral parameters. $[[R,S,T]]
= 0$ then denotes the coloured or spectral-parameter dependent QYBE
\beq \label{cqybe} R_{12}(u,v) \ R_{13}(u,w) \ R_{23}(v,w) \  = \
R_{23}(v,w) \
 R_{13}(u,w) \  R_{12}(u,v)\eeq
Equation (\ref{cqybe}) is also known as the two-parameter form of
the QYBE since the matrix $R$ depends upon two spectral parameters,
and reduces to the one-parameter form \beq R_{12}(u) \  R_{13}(u+v)
\ R_{23}(v) \ = R_{23}(v) \ R_{13}(u+v) \  R_{12}(u) \eeq for
$R(u,v)=R(u-v)$ and to the constant QYBE (\ref{qybe}) for $R(u,v)=R$

Let us now have a look at a Yang--Baxter system. A system of linear
maps
$$W: V\ot V\ \rightarrow V\ot V,\quad Z: V'\ot V'\ \rightarrow V'\ot
V',\quad X: V\ot V'\ \rightarrow V\ot V',$$ is called a
$WXZ$--system or a {\em Yang--Baxter system} \cite{h2} if the
following conditions hold: \beq \label{ybsdoub} [W,W,W] = 0 \qquad
[Z,Z,Z] = 0 \qquad [W,X,X] = 0 \qquad [X,X,Z] = 0\eeq In \cite{v} it
was observed that $WXZ$--systems with invertible $W,X$ and $Z$ can
be used to construct dually paired bialgebras of the FRT type
leading to quantum doubles. The above is one type of a constant
Yang--Baxter system that has recently been studied in \cite{np} and
also shown to be closely related to entwining structures \cite{bn}.
Other types of constant Yang--Baxter systems are those that are
related to quantised braid groups \cite{h1} and the generalised
reflection algebras \cite{c}. A number of solutions are known for
the first two types but none for the third type. Yet another type of
Yang--Baxter system is the 'coloured' Yang--Baxter system \cite{h3}
for which hardly any solutions are known.

In this paper, we provide new construction of solutions (with
Yang--Baxter operators coming from algebra structures) for
Yang--Baxter systems related to the generalised reflection algebras
as well as the coloured ones.

\bigskip

\section{Constant Yang--Baxter systems for generalised reflection algebras}

In the quantisation of nonultralocal models \cite{fm}, a quantum
version of the relations for the monodromy matrices $T$ was given in
the form
$$ A_{12} \ T_1 \ B_{12} \ T_2 \ = \ T_2 \ C_{12} \ T_1 \ D_{12} $$
where $A,B,C,D$ are numerical $n^2 \times n^2$ matrices and
$T_1=T\ot\bf{1}$, $T_2={\bf 1}\ot T$. This is a reflection-type
algebra \cite{c} that has been used for the description of open spin
chains \cite{s}. In this spirit, we consider the algebra (introduced
in \cite{fm}) generated by elements $L^k_j$, $j,k \in \{1,2,...,n\}$
satisfying the relations \beq \label{alblcons} A_{12} \ L_1 \ B_{12}
\ L_2 \ = \ L_2 \ C_{12} \ L_1 \ D_{12} \eeq where $ L \ = \ {
\{L^k_j \} }^n_{j,k=1} $. These types of algebras include well-known
algebras such as the reflection algebras, quantised function
algebras and braid groups among others. The algebra (\ref{alblcons})
has to satisfy certain consistency conditions given in terms of the
following Yang--Baxter system: \beq \label{ybsrefl}
\begin{array}{lll}
& [A,A,A]=0, \qquad & [D,D,D]=0\\
& [A,C,C]=0, \qquad & [D,B,B]=0\\
& [A, B^+, B^+ ]=0, \qquad & [D, C^+ , C^+]=0\\
& [A, C , B^+ ]=0, \qquad & [D, B , C^+]=0
\end{array}
\eeq
where $X^+ = PXP$, and $P$ is the permutation matrix.\\

The next theorem presents new solutions for this Yang--Baxter
system.

\begin{theorem}
Let $X$ be a commutative k-algebra and $ \lambda, \lambda' \in k$.

The following is a Yang--Baxter system:

$ \ A,\  B, \ C, \ D: X \ot X \rightarrow X \ot X, $

$ \ A( a \ot b)= \lambda 1 \ot ab + ab \ot 1 - b \ot a $,

$ \ B( a \ot b)= C(a\ot b)=  1 \ot ab + ab \ot 1 - b \ot a $,

$ \ D( a \ot b)= \lambda' 1 \ot ab + ab \ot 1 - b \ot a $.

\end{theorem}

\begin{proof}
$ [A,A,A]=0 $ and $ [D,D,D]=0 $ follow from Theorem 1.1 of
\cite{dn}.

$ [A,C,C]=0 $ and $ [D,B,B]=0 $ follow from Theorem 5.2 of \cite{np}.

Notice that
$B^+ (a \ot b) = PB(b \ot a)=P(1 \ot ba + ba \ot 1 - a \ot b)ba \ot 1 + 1 \ot ba - b \ot a = ab \ot 1 + 1 \ot ab - b \ot a = B(a \ot b)$.

Other relations follow immediately.
\end{proof}

Next, we look at the explicit form of the above solution in
dimension two. We consider the algebra $ \ca = \frac{
k[X]}{(X^2-\sigma)} $, where $ \sigma \in \{0, 1\}$ is a scalar.
Then $\ca$ has the basis $\{1, x \}$, where $x$ is the image of $X$
in the factor ring. We consider the basis $ \{ 1 \otimes 1, 1
\otimes x, x \otimes 1, x \otimes x \} $ of $ \ca \otimes \ca $ and
represent the operator $A$ of Theorem 2.1 in this basis:
\begin{eqnarray*}
&& A(1\otimes 1) = \lambda 1 \otimes 1 \\
&& A(1\otimes x) = \lambda 1 \otimes x \\
&& A(x\otimes 1) = (\lambda-1) 1 \otimes x + x \otimes 1 \\
&& A(x\otimes x) = \sigma (\lambda+1) 1 \otimes 1 - x \otimes x
\end{eqnarray*}
which in matrix form reads
\begin{equation}
A_{12}= \begin{pmatrix}
\lambda & 0 & 0 & \sigma (\lambda+1)\\
0 & \lambda & (\lambda-1) & 0\\
0 & 0 & 1 & 0\\
0 & 0 & 0 & -1
\end{pmatrix}
\end{equation}
Similarly,
\begin{equation}
B_{12} = C_{12} = \begin{pmatrix}
1 & 0 & 0 & 2\sigma\\
0 & 1 & 0 & 0\\
0 & 0 & 1 & 0\\
0 & 0 & 0 & -1
\end{pmatrix}
\end{equation}
\begin{equation}
D_{12}= \begin{pmatrix}
\lambda' & 0 & 0 & \sigma (\lambda'+1)\\
0 & \lambda' & (\lambda'-1) & 0\\
0 & 0 & 1 & 0\\
0 & 0 & 0 & -1
\end{pmatrix}
\end{equation}
Note that $A$ and $D$ satisfy the constant QYBE (\ref{qybe}). The
generators $L^k_j$ of the algebra (\ref{alblcons}) can be arranged
in the $2\times 2$ matrix
$$L = \begin{pmatrix}a & b\\ c & d\end{pmatrix}$$
and $L_1=L\ot\bf{1}$, $L_2={\bf 1}\ot L$ are now $4\times 4$
matrices. Then, the commutation relations among the generators
$\{a,b,c,d\}$ are obtained by evaluating both sides of the matrix
equation (\ref{alblcons}). This yields two cases (with $\lambda\neq
-1,\lambda'\neq -1$): either $\lambda = \lambda'$ or $ac=0$. In the
case $\lambda=\lambda'$, we obtain \beq \label{case1}
\begin{array}{ll}
& a^2 = 0, \quad c^2 = 0 \\
& ca = \lambda ac,\quad cb = -\lambda bc, \quad cd = \lambda dc\\
& [a,b]=2\sigma ac - (\lambda + 1)\sigma dc, \quad [a,d] = (\lambda - 1)bc\\
& [d,b] = (\lambda + 1)\sigma ac - 2\sigma cd, \quad 2\sigma ad +
b^2 - \sigma d^2 = 0
\end{array}
\eeq For $ac=0$ case, the commutation relations are \beq
\label{case2}
\begin{array}{ll}
& ac=0=ca, \quad a^2=0=c^2\\
& cb=-\lambda' bc, \quad cd=\lambda' dc\\
& ab=-\lambda' ba, \quad ad=\lambda' da\\
& [d,b] = - 2\sigma cd \\
& da=bc, \quad ba=\sigma dc\\
& 2\sigma ad + b^2 - \sigma d^2 = 0
\end{array}
\eeq It is curious to note that if $ac\neq 0$ and $\lambda+1\neq 0$,
$\lambda'+1\neq 0$, then the equality (\ref{alblcons}) requires that
$\lambda$ be equal to $\lambda'$ thus leading to the independent
relations (\ref{case1}). This is consistent with Theorem 2.1 since
if $\lambda=\lambda'$ then the operator $A=D$ and we still obtain a
solution to the Yang--Baxter system (\ref{ybsrefl}) as a special
case. On the other hand, if $\lambda\neq \lambda'$ and
$\lambda+1\neq 0$, $\lambda'+1\neq 0$, then (\ref{alblcons})
requires that $ac=ca=0$ leading to relations (\ref{case2}). It is
intriguing that, in this case, $\lambda$ disappears from the
computations and the algebra depends solely on $\lambda'$ and
$\sigma$.

\section{Coloured Yang--Baxter systems}

The commutation relations for the elements of the quantised Lax
operators depending on spectral parameters are given by \beq
A_{12}(u,v) \ L_1(u) \ B_{12}(u,v) \ L_2(v) \ = \ L_2(v) \
C_{12}(u,v) \ L_1(u) \ D_{12}(u,v) \eeq These algebras were given in
\cite{fm} and were proved useful in finding a commuting subalgebra
that facilitates the construction of quantum Hamiltonian of a model.
Here $A,B,C,D$ are spectral-parameter dependent $n^2 \times n^2$
matrices satisfying the coloured Yang--Baxter system \cite{h3} \beq
\begin{array}{lll}
& [[A,A,A]]=0, \qquad & [[D,D,D]]=0\\
& [[A,C,C]]=0, \qquad & [[D,B,B]]=0\\
& [[A, B^{++}, B^{++} ]]=0, \qquad & [[D, C^{++} , C^{++}]]=0\\
& [[A, C , B^{++} ]]=0, \qquad & [[D, B , C^{++}]]=0
\end{array}
\eeq where $X^{++}(u,v):= PX(v,u)P$. We present here a new method to
obtain families of solutions for the above system, starting from
algebra structures.\\

Let $X$ be a commutative $k$-algebra. We consider the maps $ A, C :
X \ot X \rightarrow X \ot X $ defined as follows:

$$ A(u,v)(a \ot b) = \alpha (u,v) 1 \ot ab + \beta (u,v) ab \ot 1 - \gamma (u,v)
b \ot a$$
$$C(u,v)(a \ot b) = \eta (u,v) 1 \ot ab + \zeta (u,v) ab
\ot 1 - \delta (u,v) b \ot a $$

where $ \alpha, \beta, \gamma, \eta, \zeta, \delta $ are $k$-valued functions
on $k \times k$.

We impose the conditions $ [[A,A,A]]=0 $ and $ [[A,C,C]]=0$, and
find functions which satisfy a system of functional equations.
According to \cite{np}, the condition $ [[A,A,A]]=0 $ implies that
the functions $ \alpha, \beta, \gamma $ satisfy the following system
of equations:
\begin{eqnarray} &&
(\beta(v,w)-\gamma(v,w))(\alpha(u,v)\beta(u,w) -
\alpha(u,w)\beta(u,v))\nonumber \\ &&\quad \quad \quad +
(\alpha(u,v)-\gamma(u,v))(\alpha(v,w)\beta(u,w) - \alpha(u,w)\beta(v,w))
= 0 \label{e1} \\ \nonumber \\ &&
\beta(v,w)(\beta(u,v)-\gamma(u,v))(\alpha(u,w)-\gamma(u,w)) \nonumber \\
&&\quad \quad \quad +
(\alpha(v,w)-\gamma(v,w))(\beta(u,w)\gamma(u,v)-\beta(u,v)\gamma(u,w)) 0 \label{e2} \\ \nonumber \\ && \alpha(u,v)
\beta(v,w)(\alpha(u,w)-\gamma(u,w)) + \alpha(v,w)\gamma(u,w)
(\gamma(u,v) - \alpha(u,v)) \nonumber \\ &&\quad \quad \quad +
\gamma(v,w) (\alpha(u,v)\gamma(u,w)-\alpha(u,w)\gamma(u,v)) = 0
\label{e3} \\ \nonumber \\ && \alpha(u,v)
\beta(v,w)(\beta(u,w)-\gamma(u,w)) + \beta(v,w)\gamma(u,w) (\gamma(u,v)
- \beta(u,v)) \nonumber \\ &&\quad \quad \quad + \gamma(v,w)
(\beta(u,v)\gamma(u,w)-\beta(u,w)\gamma(u,v)) = 0 \label{e4} \\
\nonumber \\ && \alpha(u,v)(\alpha(v,w)-\gamma(v,w))(\beta(u,w) -
\gamma(u,w)) \nonumber \\ &&\quad \quad \quad + (\beta(u,v)-
\gamma(u,v))( \alpha(u,w) \gamma(v,w) - \alpha(v,w) \gamma(u,w)) = 0 \
\label{e5} \end{eqnarray}

The condition $ [[A,C,C]]=0 $ implies that the functions $ \alpha,
\beta, \gamma, \eta, \zeta, \delta  $ satisfy the following system
of equations:

\begin{eqnarray} &&
(\alpha(u,v)-\gamma(u,v))(\zeta(u,w)\eta(v,w) -
\eta(u,w)\zeta(v,w))\nonumber \\ &&\quad \quad \quad +
(\alpha(u,v) \zeta(u,w) - \beta(u,v) \eta(u,w))(\zeta(v,w) - \delta(v,w))
= 0 \label{e6} \\ \nonumber \\ &&
(\gamma(u,v)-\beta(u,v))(\eta(u,w) \delta(v,w) - \delta(u,w) \eta(v,w)) \nonumber \\
&&\quad \quad \quad +
\alpha(u,v) (\delta(u,w) - \zeta(u,w))( \eta(v,w)- \delta(v,w)0 \label{e7} \end{eqnarray}

\begin{eqnarray} &&
\gamma(u,v)(\zeta(u,w) \delta(v,w) - \delta(u,w) \zeta(v,w))
+ \beta(u,v)\delta(u,w)
(\zeta(v,w) - \delta(v,w)) \nonumber \\ &&\quad \quad \quad +
\alpha(u,v) (\delta(u,w) - \zeta(u,w)) \zeta(v,w)= 0
\label{e8} \end{eqnarray}

\begin{eqnarray} &&
(\beta(u,v)-\gamma(u,v)) (\eta(u,w) - \delta(u,w)) \zeta(v,w)
\nonumber \\ &&\quad \quad \quad +
(\beta(u,v) \delta(u,w) - \gamma(u,v) \zeta(u,w))(\delta(v,w)- \eta(v,w))
= 0 \label{e9} \end{eqnarray}

\begin{eqnarray} &&
\gamma(u,v)(\delta(u,w) \eta(v,w) - \eta(u,w) \delta(v,w)) +
\alpha(u,v) \delta(u,w)(\delta(v,w) - \eta(v,w))
\nonumber \\ &&\quad \quad \quad +
\alpha(u,v) ( \eta(u,w) - \delta(u,w)) \zeta(v,w)=0
\label{e10} \end{eqnarray}

Using Theorem 2.1 from \cite{np} and by simplifying the computations,
we obtain the following solutions for the system of equations
(\ref{e1}) -- (\ref{e10}):

1)  $ \alpha(u,v)= p(u-v), \ \beta(u,v)= q(u-v), \ \gamma(u,v)= pu - qv, \
\eta(u,v)= pu - q'v, \ \zeta(u,v)= qu - p'v, \ \delta(u,v) = pu - p'v $, where
$ p, p', q, q' \in k $;

2) $ \alpha(u,v)= p(u-v), \ \beta(u,v)= q(u-v), \ \gamma(u,v)= pu - qv, \
\eta(u,v)= p(\lambda u - \mu v), \ \zeta(u,v)= q(\lambda u - \mu v)  , \
\delta(u,v) = p \lambda u - q \mu v $, where
$ p, q, \lambda, \mu \in k $.\\

These solutions lead to the following theorem:

\begin{theorem}
Let $X$ be a commutative k-algebra and $ p, p', q, q', \lambda, \mu
\in k$.

The coloured operators $ \ A,\  B, \ C, \ D: k \ot k \rightarrow
End_k ( X \ot X) $ are solutions for the coloured Yang--Baxter
system in the following
two cases:\\

1)
$ \ A(u,v)( a \ot b)= p(u-v) 1 \ot ab + q(u-v) ab \ot 1 - (pu-qv) b \ot a $,

$ \ B(u,v)(a \ot b)= (p'u-qv) 1 \ot ab + (q'u-pv) ab \ot 1 - (p'u-pv) b \ot a $,

$ \ C(u,v)(a\ot b)= (pu-q'v) 1 \ot ab + (qu-p'v) ab \ot 1 - (pu-p'v) b \ot a $,

$ \ D(u,v)( a \ot b)= p'(u-v) 1 \ot ab + q'(u-v) ab \ot 1 - (p'u-q'v) b \ot a $;\\

2)
$ \ A(u,v)( a \ot b)= p(u-v) 1 \ot ab + q(u-v) ab \ot 1 - (pu-qv) b \ot a $,

$ \ B(u,v)(a \ot b)= q(\mu u - \lambda v) 1 \ot ab +
p(\mu u - \lambda v) ab \ot 1 - (q \mu u - p \lambda v) b \ot a $,

$ \ C(u,v)(a\ot b)=  p(\lambda u - \mu v) 1 \ot ab + q(\lambda u - \mu v) ab \ot 1 - (p \lambda u - q \mu v)  b \ot a $,

$ \ D(u,v)( a \ot b)= q(u-v) 1 \ot ab + p(u-v) ab \ot 1 - (qu-pv) b \ot a $.

\end{theorem}

\begin{proof}
Let us observe that $ C^{++}(u,v) = - B(u,v) $ in both cases of the
theorem. Everything now follows from the above analysis.
\end{proof}

Finding other solutions for the system of equations (\ref{e1}) -- (\ref{e10}) is an open problem.
We now present the above solutions in dimension two. Consider the
algebra $\ca$ of the previous section and working in the same basis
for $\ca\ot \ca$, we obtain the following matrix solutions for case
1 of Theorem 3.1: \beq
A(u,v)= \begin{pmatrix}
qu-pv & 0 & 0 & \sigma (q+p)(u-v)\\
0 & p(u-v) & (q-p)v & 0\\
0 & (q-p)u & q(u-v) & 0\\
0 & 0 & 0 & qv-pu
\end{pmatrix}
\eeq \beq B(u,v)= \begin{pmatrix}
q'u-qv & 0 & 0 & \sigma [(p'+q')u-(p+q)v]\\
0 & p'u-qv & (p-q)v & 0\\
0 & (q'-p')u & q'u-pv) & 0\\
0 & 0 & 0 & pv-p'u
\end{pmatrix}
\eeq \beq C(u,v)= \begin{pmatrix}
qu-q'v & 0 & 0 & \sigma [(p+q)u-(p'+q')v]\\
0 & pu-q'v & (p'-q')v & 0\\
0 & (q-p)u & qu-p'v) & 0\\
0 & 0 & 0 & p'v-pu
\end{pmatrix}
\eeq \beq D(u,v)= \begin{pmatrix}
q'u-p'v & 0 & 0 & \sigma (q'+p')(u-v)\\
0 & p'(u-v) & (q'-p')v & 0\\
0 & (q'-p')u & q'(u-v) & 0\\
0 & 0 & 0 & q'v-p'u
\end{pmatrix}
\eeq Similarly, for case 2 we have \beq A(u,v)= \begin{pmatrix}
qu-pv & 0 & 0 & \sigma (q+p)(u-v)\\
0 & p(u-v) & (q-p)v & 0\\
0 & (q-p)u & q(u-v) & 0\\
0 & 0 & 0 & qv-pu
\end{pmatrix}
\eeq \beq B(u,v)= \begin{pmatrix}
p\mu u-q\lambda v & 0 & 0 & \sigma (p+q)(\mu u-\lambda v)\\
0 & q(\mu u-\lambda v) & (p-q)\lambda v & 0\\
0 & (p-q)\mu u & p(\mu u-\lambda v) & 0\\
0 & 0 & 0 & p\lambda v-q\mu u
\end{pmatrix}
\eeq \beq C(u,v)= \begin{pmatrix}
q\lambda u-p\mu v & 0 & 0 & \sigma (p+q)(\lambda u-\mu v)\\
0 & p(\lambda u-\mu v) & (q-p)\mu v & 0\\
0 & (q-p)\lambda u & q(\lambda u-\mu v) & 0\\
0 & 0 & 0 & q\mu v-p\lambda u
\end{pmatrix}
\eeq \beq D(u,v)= \begin{pmatrix}
pu-qv & 0 & 0 & \sigma (q+p)(u-v)\\
0 & q(u-v) & (p-q)v & 0\\
0 & (p-q)u & p(u-v) & 0\\
0 & 0 & 0 & pv-qu
\end{pmatrix}
\eeq

\section{Conclusions}

In this work, we have investigated new constructions of the constant
and coloured Yang--Baxter systems from the viewpoint of Yang--Baxter
operators from algebra structures. A solution for the constant
Yang--Baxter system related to the generalised reflection algebra is
presented and the commutation algebra structure in dimension two is
exhibited in detail. Furthermore, for the coloured Yang--Baxter
systems, we have given a construction that involves solving a
nontrivial system of functional equations. We obtain two families of
such solutions and it remains an open problem to find and classify
solutions associated to a certain system of functional equations.

\end{document}